\documentclass[a4paper,12pt,twoside,tbtags]{article}
\usepackage{amstext,amsmath,amssymb,mathrsfs}
\usepackage{amsthm}
\usepackage{bm}
\usepackage[dvipdfmx]{graphicx}

\def\bbC{\mathbb{O}}%\def\bbC{\mathbb{C}}
\def\D{\mathbb{D}}

\def\I{\mathbb{I}}

\def\M{\mathbb{M}}
\def\N{\mathbb{N}}
\def\P{{\mathbb P}}

\def\R{\mathbb{R}}

\def\Z{\mathbb{Z}}

\def\C{O}%\def\C{\mathbf{C}}

\def\bP{{\bf P}}

\def\1{{\bf 1}}

\def\b{{\bm b}}

\def\u{{\bm u}}

\def\w{{\bm w}}
\def\x{{\bm x}}
\def\y{{\bm y}}

\def\B{{\bm B}}
\def\L{{\bm L}}
\def\H{{\bm H}}
\def\S{\bm{S}}
\def\X{{\bm X}}
\def\Y{{\bm Y}}

\def\cB{\mathcal{B}}
\def\cF{\mathcal{F}}

\def\ml{\mathfrak{l}}
\def\mfu{\mathfrak{u}}
\def\mH{\mathfrak{H}}
\def\mI{\mathfrak{I}}
\def\mM{\mathfrak{M}}

\def\mX{\mathfrak{X}}

\def\={\stackrel{\rm (law)}{=}}

\newcommand{\As[1]}{\textbf{(#1)}}

\newtheorem{thm}{Theorem}[section]

\newtheorem{lmm}[thm]{Lemma}

\newtheorem{prp}[thm]{Proposition}

\newtheorem{dfn}[thm]{Definition}
\newtheorem{exa}[thm]{Example}
\newtheorem{rem}{Remark}%%%

\newcommand{\SSC}[1]{\section{#1}\setcounter{equation}{0}}
\newcommand{\SSSC}[1]{\subsection{#1}\setcounter{equation}{0}}

\begin{document}
%\mainmatter
\title{\bfseries Infinite particle systems with hard-core \\ and long-range interaction}
%\titlerunning{Infinite particle systems with hard-core and long-range interaction}

\author{Hideki Tanemura\thanks{Department of Mathematics, Keio University, Yokohama 223-8522, Japan.
e-mail: tanemura@math.keio.ac.jp}}
%\authorrunning{Hideki Tanemura}
%\tocauthor{David Berger, Franziska K\"uhn and Ren\'e L. Schilling}
%\institute{ Keio University \\
%\email{tanemura@math.keio.ac.jp}  }
 
\maketitle

\begin{abstract}  
A system of Brownian hard balls is regarded as a reflecting Brownian motion in the configuration space 
and can be represented by a solution to a Skorohod-type equation.
In this article, we consider the case that there are an infinite number of balls, and the interaction between balls is given by the long-range pair interaction. We discuss the existence and uniqueness of strong solutions to the infinite-dimensional Skorohod equation.

%\keywords{Keywords} Stochastic differential equations, Infinite particle systems, Skorohod equations
%\subjclass{60K35, 60J46, 60J60    }

%\ack{ The author is supported in part by JSPS KAKENHI Grant Numbers JP.16H06338, JP.19H01793.}

\end{abstract}

%\bigskip\noindent

%%%%%%%%%%%%%%%%%%%%%%%
%\section{Introduction}
\SSC{Introduction}
%%%%%%%%%%%%%%%%%%%%%%%%

In this article, we study systems of interacting Brownian motions on $\R^d$, $d\ge 2$. 
Let $\Phi : \R^d \to (-\infty, \infty]$
be a self (free) potential and
$\Psi : \R^d\times \R^d \to (-\infty, \infty)$
be a symmetric pair-interaction potential.
In the case that these potentials are smooth,
the system is described by the stochastic differential equation (SDE)
$$
dX_t^j = dB_t^j -\frac{1}{2}\nabla\Phi(X_t^j)dt -\frac{1}{2}\sum_{k\in \Lambda, k\not=j} \nabla\Psi (X_t^j,X_t^k)dt, \quad j\in\I,
$$
where $B_t^j$, $j\in\Lambda$ are independent Brownian motions and $\I$ is a countable index set. 
We consider in this article the case that $\Phi$ is smooth and $\Psi$ is a pair potential with a hard core of radius $r>0$ (i.e., $\Psi=\Psi_{\rm hard} + \Psi_{\rm sm}$):
\begin{align}
&\Psi_{\rm hard}(x,y)=
\begin{cases}0 \quad \mbox{ if \; $|x-y|\ge r$,}
\\
\infty  \quad \mbox{if \; $|x-y|< r$,}
\end{cases}
\mbox{ the hard-core pair potential, }
\nonumber
\\
&\Psi_{\rm sm}(x,y)=\Psi_{\rm sm}(x-y) : \mbox{a translation invariant smooth potential.}
\nonumber
\end{align}
The system can be regarded as that of balls with a radius $r>0$.

When $\I$ is a finite subset of $\N$,
Saisho\cite{Saisho} and Saisho and Tanaka\cite{Saisho-Tanaka} showed that a system of interacting Brownian balls with potential $\Phi$ and $\Psi_{\rm sm}$ can be represented by the unique solution of the Skorohod-type equation
\begin{align}
&dX_t^j=dB_t^j -\frac{1}{2}\nabla\Phi(X_t^j)dt -\frac{1}{2}
\sum_{k\in\I, k\not= j}\nabla \Psi_{\rm sm}(X_t^j-X_t^k)dt
\nonumber
\\
&\qquad \qquad+ \sum_{k\in\I, k\not= j}(X_t^j-X_t^k)dL_t^{jk},
 \quad j\in\I,
\tag{SKE-$\I$}\label{SKE_I}
\\ \nonumber 
&|X_t^j -X_t^k| \ge r, \quad j,k \in\I,
\end{align}
where $L_t^{jk}$ $j,k\in\I$ are non-decreasing functions satisfying
\begin{equation}\nonumber
L_t^{jk} = L_t^{kj}= \int_0^t \mathbf{1}(|X_s^j-X_s^k|=r)dL_s^{jk},
 \quad j,k\in\I.
\end{equation}
For $\I=\N$, the existence and uniqueness of solutions of (SKE-$\N$) have been solved in the cases that $\Phi=\Psi_{\rm sm}=0$ \cite{T96} and $\Phi=0$  and $\Psi_{\rm sm}$ has stretched exponential decay \cite{FRT}. 
In these cases, the interaction among particles has a short range.
%These results were obtained under the condition that the interaction among particles is short range.
The purpose of this article is to generalize the results for long-range interaction including the case that $\Psi_{\rm sm}$ has polynomial decay.

Let $\mX$ be the configuration space of unlabeled balls with diameter $r>0$. 
The space $\mX$ is a compact Polish space with the vague topology. 
Using Dirichlet form theory, we can construct an $\mX$-valued process $\Xi$ describing an interacting system with an infinite number of unlabeled particles \cite{O96,o.rm} including the case with a hard-core interaction. See, for example, \cite{F1,F2} for the relation between Dirichlet forms and reflecting Brownian motions.
For a system with a finite number of balls, the existence of the solution to the Skorohod equation is derived through Fukushima decomposition from the process constructed using a Dirichlet form. (See, for instance, \cite{Ch93}.)  
For a system with an infinite number of balls, we need to label the balls because the coordinate function is not locally in the domain of the Dirichlet form.
To label the balls in the system, we use a sequence of tagged particle processes $\{(\X^m, \Xi^m)\}_{M\in \N}$ with consistency, as introduced by Osada\cite{o.tp}. 
Additionally, we can apply the argument in \cite{o.isde} to the case with a hard-core interaction.
We can then apply the Fukushima decomposition to our case and show the existence of a solution $\X$ for SKE-$\N$. 

The existence and uniqueness of strong solutions have been discussed \cite{o-t.tail}.
In the cited paper, we introduced an infinite system of finite-dimensional SDEs associated with a solution $\X$ of (SKE-$\N$). 
We showed that under the condition that each finite-dimensional SDE has a unique strong solution, referenced as \As[IFC], there exists a strong solution and a unique strong solution of (SKE-$\N$) under some constraints.  
However, we are not sure if the IFC holds in our model with a hard-core interaction.
In the present paper, we introduce two conditions, namely $\mI$-IFC and the finite cluster property \As[FCP].
We present a result for the existence and uniqueness of strong solutions of (SKE-$\N$) under \As[$\mI$-IFC] and \As[FCP] and verify these conditions in our setting.

Section 2 prepares notations and cites results on the construction of the unlabeled process $\Xi$ and
labeled process $\X$
and the Skorohod equation.
Section 3 presents results, Proposition \ref{prp:31} and Theorem \ref{th:main}. 
Section 4 presents the proof of Theorem \ref{th:main}.

%%%%%%%%%%%%%%%%%%%%%%%%%%%%%%%%%
%\section{Preliminaries}
\SSC{Preliminaries}
%%%%%%%%%%%%%%%%%%%%%%%%%%%%%%%%%

%%%%%%%%%%%%%%%%%%%%%%%%%%%%%%%%%%%%%%%%%%
%\subsection{Systems of unlabeled hard balls}
\SSSC{Systems of unlabeled hard balls}
%%%%%%%%%%%%%%%%%%%%%%%%%%%%%%%%%%%%%%%%%%

We denote the configuration space of unlabeled balls with radius $r>0$ in $\R^d$ by
$$
\mX=\{ \xi= \{x^j\}_{j\in\I } : |x^j-x^k|\ge r \quad j \not=k, \; \I\mbox{ is countable} \}.
$$
We can regard an element $\xi=\{x^j\}_{j\in\I}\in \mX$ as a Radon measure $\sum_{j\in\I}\delta_{x^j}$ and $\mX$ as a subset of the set $\mM$ of non-negative Radon measures:
$$
\mM =\mM(\R^d)= \Big\{ \xi(\cdot)=\displaystyle{\sum_{j\in \I}} \delta_{x^j}(\cdot):
\xi (K) < \infty, \mbox{ $\forall K\subset\R^d$ compact} \Big\}.
$$
where $\delta_x$ is the delta measure at $x$.
We remark that $\mX$ is compact with the vague topology. 
We denote by $\pi_A(\xi)$ the restriction of $\xi$ on $A\in\R^d$.
\vskip 3mm

%The index set $\I$ is countable.
%$\mathfrak{M}$ is a Polish space with the vague topology.
%$$
%\mbox{$\mN\subset \mathfrak{M}$ is  {relative compact}}
%\Leftrightarrow
%\sup_{\xi \in \mN}\xi (K) <\infty,
%\; \mbox{$\forall K\subset\R^d$ compact}
%$$

We cite the definition of the quasi-Gibbs measure on $\mM$ \cite[Definition 2.1]{o.rm2}. 
For $\zeta \in \mM$,
the Hamiltonian of $\Phi$, $\Psi$ on $U_\ell=\{x\in \R^d: |x|\le \ell \}$ is given by
$$
H_\ell(\zeta)=\sum_{x\in { \rm supp }\zeta\cap{U_\ell}}\Phi(x)+\sum_{x,y\in { \rm supp }\zeta\cap{U_\ell},x\not=y}\Psi (x,y).
$$
Let $\Lambda$ be the Poisson random field on $\R^d$ with intensity measure $dx$.
\begin{dfn}[Quasi-Gibbs state]
A probability measure $\mu$ is called a $(\Phi,\Psi)$-{quasi-Gibbs state}
if %the regular conditional distribution 
$$
\mu_{\ell,\xi}^m(d\zeta)=\mu(\pi_{U_\ell}(\zeta)\in d\zeta | \pi_{U_\ell^c}(\xi)=\pi_{U_\ell^c}(\zeta),\zeta(U_\ell)=m)
$$
satisfies that for $\ell,m,k\in\N$, $\mu$-a.s. $\xi$, 
$$
c^{-1}e^{-H_\ell(\zeta)}\Lambda_\ell^m(d\zeta) \le 
\mu_{\ell,\xi}^m(d\zeta)
\le c e^{-H_\ell(\zeta)}\Lambda_\ell^m(d\zeta), 
\text{}
$$
where 
$c=c(\ell,m,\xi)>0$ is a constant depending on $\ell,m,\xi$ and
$$
\Lambda_\ell^m (\cdot) = \Lambda( \pi_{U_\ell} \in \cdot|\mathfrak{M}_\ell^m) \quad \text{with} \quad
\mathfrak{M}_\ell^m=\{\xi(U_\ell)=m\}.
$$
\end{dfn}

A function $f$ on $\mX$ is called a {polynomial function} if it can be expressed as
$$
f(\xi) =Q \left( \langle \varphi_1,\xi\rangle, \langle \varphi_2,\xi\rangle, \dots, \langle \varphi_\ell,\xi\rangle\right) 
$$
with a polynomial function $Q$ on $\R^\ell$ and smooth functions $\varphi_j$, $1\le j\le \ell$,
on $\R^d$ with compact support, where
$$
\langle \varphi, \xi\rangle = \int_{\R^d}\varphi(x) \xi(dx).
$$ 
We denote by ${\cal P}$ the set of all polynomial functions on $\mM$.
A polynomial function is a local and smooth function; i.e., there is a compact set $K$ such that
$$
f(\xi)=f(\pi_K(\xi))
$$
and symmetric smooth functions $\hat{f}_n$, $n\in\N$ such that
$$ 
f(\xi) = \hat{f}(x_1,\dots, x_n), \qquad
\mbox{if $\xi\cap K= \sum_{j=1}^n \delta_{x^j}$}.
$$
The sequence of the functions $\hat{f}_n$, $n\in\N$ is called a K-representation of $f$.

For $f\in \mathcal{P}$, we introduce the {square field} on $\mM$ defined by
\begin{align}\label{:21d}
\mathbb{D}(f,g)(\xi)=\frac{1}{2}\int_{\R^d}\xi(dx)
\nabla_{x} f(\xi) \cdot 
\nabla_{x} g(\xi).
\end{align}
For a probability measure $\mu$ on $\mX$, we introduce the {bilinear form} on $L^2(\mu)$ defined by
\begin{align}\nonumber
&{\cal E}^{\mu}(f,g)= \int_{\mathfrak{M}}
\mathbb{D}(f,g)(\xi)
\mu(d\xi),
\quad f,g \in \mathcal{D}_{\circ}^{\mu},
\\ \nonumber
&\mathcal{D}_{\circ}^{\mu}=\{ f\in \mathcal{P} : \parallel f \parallel_1 <\infty \},
\end{align}
where
$$
\parallel f \parallel_1^2=\parallel f\parallel_{L^2(\mu)}^2+{\cal E}^{\mu}(f,f). 
$$

We make the following assumptions.

\vskip 3mm

\noindent \As[A0] The pair potential $\Psi$ has a hard-core interaction with radius $r>0$; i.e.,
$\Psi (x,y) = \infty$ if $|x-y|\le r$.

\vskip 3mm

\noindent \As[A1] $\mu$ is a $(\Phi,\Psi)$-quasi-Gibbs state, and 
$\Phi:\R^d\to\R\cup \{\infty\}$ and $\Psi:\R^d\times\R^d\to\R\cup \{\infty\}$ satisfy
\begin{eqnarray}
&&c^{-1}\Phi_0(x)
\le \Phi(x)
\le c \; \Phi_0(x),
\nonumber\\
&&c^{-1}\Psi_0(x-y)
\le \Psi(x,y)
\le c \; \Psi_0(x-y)
\nonumber
\end{eqnarray}
for some $c>1$ and are locally bounded from below and upper semi-continuous functions $\Phi_0, \Psi_0$ with $\{x\in\R^d : \Psi_0(x) =\infty \}$ being compact. 

From \As[A0], $\mu$ satisfies $\mu(\mX)=1$.
We cite the result in \cite[Lemma 2.1]{o.rm} with \cite{o-t.core}.

\begin{lmm}[\cite{o.rm,o-t.core}]
Assume \As[A0] and \As[A1].
\\ \noindent
{\rm (i)}  $(\mathcal{E}^\mu, \mathcal{D}_{\circ}^\mu)$ is closable on $L^2(\mX, \mu)$, 
and the closure $(\mathcal{E}^\mu, \mathcal{D}^\mu)$ is a regular Dirichlet form. 
\\ \noindent
{\rm (ii)} There is a diffusion process $(\Xi_t, \mathbb{P} _{\xi})$ associated with $(\mathcal{E}^\mu, \mathcal{D}^\mu)$ on $L^2(\mu)$. 
\\ \noindent
{\rm (iii)} $(\Xi_t, \mathbb{P} _\mu)$ is a reversible process, where $\P_{\mu}=\int_{\mX}\mu(d\xi) \P_{\xi}$.
\end{lmm}

\vskip 3mm

%\begin{rem} Suppose that
%Gibbs states of Ruelle class are quasi-Gibbs state. Then
%for Lennard-Jones 6-12 potential :
%$\Psi_{\rm sm}=\Psi_{6,12}(x) = \{|x|^{-12}-|x|^{-6}\}.$ and Riesz potential : 
%$\Psi_{\rm sm}=\Psi_a(x)=(\beta / a )|x|^{-a}.$, $ d < a \in \N  $ we can apply the above theorem.
%\end{rem}

%%%%%%%%%%%%%%%%%%%%%%%%%%%%
\subsection{Systems of labeled balls}

We denote the unlabeled configuration space of an infinite number of balls by
$$
\mX_{\infty}=\{ \xi= \{x^j\}_{j\in\N } : |x^j-x^k|\ge r \quad j \not=k \}.
$$
We make the following assumption because we are studying a system of infinite particles.

\vskip 3mm

\noindent \As[A2] $\mu(\mX_\infty)=1$. 

\vskip 3mm

We denote the configuration space of labeled balls by
\begin{align}\notag%\label{:22a}
\S_{\rm hard}=\{
\x=(x^j)_{j\in\N} \in (\R^d)^{\N} : |x^j -x^k|\ge r, \quad j\not= k 
\}.
\end{align}
We introduce the {unlabel} map $\mfu: \S_{\rm hard} \to \mX_{\infty}$ defined by
\begin{align}\label{:22c}
\mfu((x^j)_{j\in\N}) =\{x^j\}_{j\in\N}
\end{align}
and a {label} map $\ml: \mX_{\infty} \to \mathbf{S}_{\rm hard}$ such that
\begin{align}\label{:22d}
\ml(\xi)=(x^j)_{j\in\N}, \quad \mbox{if $\xi=\{x^j\}_{j\in\N}\in \mX_\infty$.} 
\end{align}

From the hard-core condition of the configuration space $\mX$ and the continuity of the trajectory of the process $\Xi$, we can lift the unlabeled dynamics $\Xi= \{ X^j\}_{j\in\N}$ to labeled dynamics $\X=(X^j)_{j\in\N}$. Here, $X^j$ is an $\R^d$-valued continuous process on one of the intervals of the form $[0, b)$ and $(a, b)$, $0< a < b\le \infty$.
We refer to $X^j$ as a tagged particle.
If $b<\infty$, we say the tagged particle explodes.
If $a>0$, we say that the tagged particle enters.
We make the following assumption.

\vskip 3mm

\noindent \As[NEE] \quad $\Xi_t$ is an $\mX$-valued diffusion process in which no tagged particle explodes
or enters.

\vskip 3mm
For a topological space $S$, $W(S)$ denotes the set of continuous paths from $\R_+ := [0,\infty)$ to $S$ and we put $W_x(S)= \{ w\in W(S) : w(0)=x\}$ for $x\in S$.
Under {\bf (NEE)}, we can construct a labeled map $\ml_{\rm path}$ from $W( \mX_\infty)$ to $W(\mathbf{S}_{\rm hard})$ such that for $\Xi = \{ X^j \}_{j\in \N} \in W( \mX_\infty)$ we have
$$
\ml_{\rm path}(\Xi)= (X^j)_{j\in\mathbb{N}}\equiv \X.
$$
We remark that $\ml_{\rm path}(\Xi)_t \not= \ml(\Xi_t)$. We also put for $m\in\N$
$$
\ml_{\rm path}^{[m]}(\Xi)= ( (X^j)_{j=1}^m, \{X_j\}_{j=m+1}^\infty) 
=: (\X^m, \Xi^m).
$$

We quote the results in \cite[Theorem 2.5]{o.tp}.
See also \cite[Lemma 1.2]{o-t.tail}.

\begin{lmm}[\cite{o.tp}]
Assume \As[A0]--\As[A2]. $(\Xi_t, \P_{\mu})$ satisfies \As[NEE].
\end{lmm}

From the above lemma we can lift the unlabeled process $(\Xi, \P_\xi)$ to a labeled process $(\X, \bP_{\x})$ such that
%$$
%\Xi_t= \mfu(\X_t),
%\quad \P_{\mfu(\x)}= \bP_{\x}\circ \mfu^{-1}
%\quad \mbox{and} \quad \ml(\xi)=\x.
%$$
$$
\X= \ml_{\rm path}(\Xi),
\quad \bP_{\x} = \P_{\mfu(\x)}
\quad \mbox{and} \quad  \x = \ml(\xi).
$$

$\ml_{\rm path}(\Xi)_t$ depends on not only $\Xi_t$ but also the trajectory of $\Xi$, 
and $\X^{m} = (X^1, X^2, \dots, X^m)$, $m\in \N$, is thus not a Dirichlet process for $\Xi$.
Then, using the argument in \cite{o.tp}, we introduce the $m$-labelled processes $(\X^{m}, \Xi^{m})$, $m\in \N \cup \{0\}$, for which $\X^{m}$ is a Dirichlet process.

We shall present the Dirichlet form associated with the $m$-labeled process. 
Let $\mu^{[m]}$ be the reduced $m$-Campbell measure on $(\R^d)^m \times \mX$ for $\mu$ defined as
\begin{align}\notag
\mu^{[m]}(d\x^m d\eta) = \rho^m (\x^m)\mu_{\x^m}(d\eta)d\x^m,
\end{align}
where $\rho^m$ is the $m$-point correlation function of $\mu$ with respect to the Lebesgue measure $d\x^m$
and $\mu_{\x^m}$ is the reduced Palm measure conditioned at $\x^m \in (\R^d)^m$. See, for instance, \cite{K17} for these definitions.
We introduce the bilinear form $(\mathcal{E}^{\mu^{[m]}}, \mathcal{D}_{\circ}^{\mu^{[m]}})$ defined by
\begin{align}\nonumber%\label{:22f}
&\mathcal{E}^{\mu^{[m]}}(f,g) = \int_{(\R^d)^m \times \mX} 
\Big\{ \frac{1}{2}\sum_{i=1}^m \frac{df}{dx^i}\frac{dg}{dx^i} + \mathbb{D}(f,g)
\Big\}
\mu^{[m]}(d\x^m d\eta),
\\ \notag
&\mathcal{D}_{\circ}^{\mu^{[m]}} = \Big\{
f\in C_0^\infty ((\R^d)^m) \otimes \mathcal{D}_\circ ; \mathcal{E}^{\mu^{[m]}}(f,f) < \infty,
f \in L^2 ( (\R^d)^m \times\mX, \mu^{[m]})
\Big\}
,\end{align}
where $\displaystyle{\frac{\partial}{\partial x^j}}$ is the nabla in $\R^d$,
and $\D$ is defined by \eqref{:21d}.

We quote the following.
\begin{lmm}[\cite{o.tp}]%%%%%%%%%%%%%%%%%%%%%%%%%%%%%%%%%%%%%%%%%%%%%%
Assume \As[A0]--\As[A2]. Let $m\in\N$.
\\
{\rm (i)} The bilinear form $(\mathcal{E}^{\mu^{[m]}}, \mathcal{D}_{\circ}^{\mu^{[m]}})$ is closable. 
Its closure, denoted by $(\mathcal{E}^{\mu^{[m]}}, \mathcal{D}^{\mu^{[m]}})$, 
is associated with the diffusion process $((\X_t^{m}, \Xi_t^{m}), \P_{(\x^m,\eta)}^{[m]})$.
\\
{\rm (ii)} The sequence $ \{ ((\X_t^{m}, \Xi_t^m), \P_{(\x^m,\eta)}^{[m]})\}_{m\in \N}$ satisfies the consistency condition
\begin{align}\notag%\label{:22g}
\P_{(\x^m,\eta)}^{[m]} = \P_{\mfu (\x^m, \eta)}\circ (\ml_{\rm path}^{[m]})^{-1},
\quad
\P_{(\x^m,\eta)}^{[m]} \circ \mfu^{-1}= \P_{\mfu (\x^m, \eta)},
\end{align}
where $\mfu (\x^m, \eta) = \{ x^j \}_{j=1}^m \cup \eta \in \mX$, if $\{ x^j \}_{j=1}^m \cap \eta = \emptyset$ and  $\{ x^j \}_{j=1}^m \cup \eta \in \mX$.
\end{lmm}%%%%%%%%%%%%%%%%%%%%%%%%%%%%%%%%%%%%%%%%%%%%%%%%%%%%%%%%%%%%

We can construct from this lemma the labeled process $\X = (X^1, X^2, \dots)$ satisfying
\begin{align}\label{:22k}
\X^{m} = (X^{1}, \dots, X^m), \quad \Xi^{m} = \{ X^j\}_{j=m+1}^\infty, \quad m\in \N
.\end{align}
In particular, we can regard the process $X^j$, $j \le m $ as a Dirichlet process of the diffusion $(\X^m, \Xi^m)$ associated with the Dirichlet form $(\mathcal{E}^{\mu^{[m]}}, \mathcal{D}^{\mu^{[m]}})$.

%%%%%%%%%%%%%%%%%%%%%%%%%%%%%%%%%%%%%%%%%%%%%%%%%%%%
\subsection{Skorohod equation}%%%%%%%%%%%%%%%%%%%%%%
%%%%%%%%%%%%%%%%%%%%%%%%%%%%%%%%%%%%%%%%%%%%%%%%%%%%%

Let $D$ be an open domain in $\R^N$, $N\in \N$.
Let $\mathcal{N}_x=\mathcal{N}_x(D)$ be the set of inward normal unit vectors at $x\in \partial D$,
$$
\mathcal{N}_x=\bigcup_{\ell>0}\mathcal{N}_{x,\ell}
\quad
\mathcal{N}_{x,\ell}
=\{\bm{n}\in\R^N: |\bm{n}|=1, U_\ell(x-\ell\bm{n})\cap D=\emptyset \}.
$$
For $x\in \overline{D}$ and $w\in W_0(\R^N)$, we consider what is called the {\it Skorohod equation}:
\begin{equation}\tag{SK}\label{Sk}
\zeta(t)=x+w(t)+\varphi(t), \quad t\ge 0.
\end{equation}
A solution to (\ref{Sk}) is a pair $(\zeta, \varphi)$ satisfying (\ref{Sk}) with the following two conditions.
\begin{itemize}
\item[(1)] $\zeta \in W(\overline{D})$.

\item[(2)] $\varphi$ is an $\R^N$-valued continuous function with bounded variation on each finite time interval satisfying $\varphi(0)=0$ and
$$
\varphi(t)=\int_0^t \bm{n}(s)d \|\varphi\|_s, \quad 
\|\varphi\|_t=\int_0^t
\mathbf{1}_{\partial D}(\zeta(s))d\|\varphi\|_s,
$$
where $\bm{n}(s)\in \mathcal{N}_{\zeta(s)}$ if $\zeta(s)\in \partial D$ and
$\|\varphi\|_t$ denotes the total variation of $\varphi$ on $[0,t]$.
\end{itemize}

We introduce the following conditions for $D$.
\begin{itemize}
\item[(A)](Uniform exterior sphere condition) There exists a constant $\alpha_0>0$ such that
$$
\forall x\in \partial D, \quad \mathcal{N}_{x}=\mathcal{N}_{x,\alpha_0}\not=\emptyset.
$$

\item[(B)] There exists constants $\delta_0>0$ and $\beta_0\in [1,\infty)$ such that for any $x\in\partial D$ there exists a unit vector $\bm{l}_x$ verifying
$$
\forall \bm{n}\in \bigcup_{y\in U_{\delta_0}(x)\cap \partial D}\mathcal{N}_{y},
\quad \langle \bm{l}_x, \bm{n}\rangle \ge \frac{1}{\beta_0}.
$$
\end{itemize}

We quote the results in \cite{Saisho,Saisho-Tanaka} and \cite[Lemmas 3.1 and 3.2]{FRT}.
%%%%%%%%%%%%%%%%%%%%%%%%%%%%%%%%%%%%%%%%%
%%%%%%%%%%%%%%%%%%%%%%%%%%%%%%%%%%%%%%%%%
\begin{lmm}[\cite{Saisho,Saisho-Tanaka,FRT}]\label{l:23a}
{\rm (i)} Suppose $D$ satisfies conditions (A) and (B).
There is then a unique solution of (\ref{Sk}).

\noindent {\rm (ii)} Suppose that $D$ satisfies conditions (A) and (B).
Let $(\zeta^{(i)}, \phi^{(i)})$ be the solution of (\ref{Sk}) for $x^{(i)}$ and $w^{(i)}$, $i=1,2.$
Then, for each $T>0$, there exists a constant $C=C(\alpha_0,\beta_0,\delta_0)$ such that
\begin{align}\label{:23c}
|\zeta^{{(1)}}(t)-\zeta^{{(2)}}(t)|\le (\|w^{(1)}- w^{(2)}\|_t +|x^{(1)}-x^{(2)}|)
e^{C(\|\varphi^{(1)}\|_t + \|\varphi^{(2)}\|_t)}
\end{align}
and 
\begin{align}\label{:23d}
\|\varphi^{(i)}\|_t\le f(\Delta_{0,T,\cdot}(w^{(i)}), \sup_{0\le s \le t}|w^{(i)}|), \quad 0\le t \le T, \; i=1,2,
\end{align}
where $f$ is a function on $W_0(\R_+)\times \R_+$ depending on $\alpha_0,\beta_0,\delta_0$,
and $\Delta_{0,T,\delta}(w)$ which denotes the modulus of continuity of $w$ in $[0,T]$.

\noindent {\rm (iii)} The configuration space of $n$ balls with diameter $r>0$,
$$
D_n=\{ \x=(x^1,x^2,\dots,x^n)\in (\R^d)^n :
|x^j -x^k| > r, \quad j\not=k \},
$$
satisfies conditions (A) and (B).
\end{lmm}
As a corollary of this lemma, 
the existence of a unique strong solution of Skorohod SDEs has been proved \cite[Theorem 5.1]{Saisho}.
We then see the existence of a unique strong solution of \eqref{SKE_I} if $\I$ is a finite subset of $\N$. 
An approximation Skorohod-type equation was introduced in the proof of \cite[Theorem 5.1]{Saisho}. 
In our setting, the equation is written as
\begin{align}
&dX_{n}^j(t)=dB_t^j -\frac{1}{2}\nabla\Phi(X_n^j(h_n(t))dt -\frac{1}{2}
\sum_{k\in\I, k\not= j}\nabla \Psi_{\rm sm}(X_n^j(h_n(t))-X_n^k(h_n(t))dt
\nonumber
\\
&\qquad \qquad+ \sum_{k\in\I, k\not= j}(X_n^j(t)-X_n^k(t))dL_n^{jk}(t),
 \quad j\in\I,
\tag{SKE(n)-$\I$}\label{SKE(n)_I}
\end{align}
with the initial condition $(X_n^j(0))_{j\in\I}=(X^j(0))_{j\in\I}$,
where 
$$
h_n(0)=0, \ h_n(t)=(k-1)2^{-n}, \ (k-1)2^{-n}< t \le k2^{-n}, \
k\in \N, \ n\in\N.
$$
By Lemma \ref{l:23a}, \eqref{SKE(n)_I} has a unique strong solution, 
and the limit  of the sequence $\{((X_n^j)_{j\in\I}, (L_n^{jk})_{j,k\in\I})\}_{n\in\N}$ as $n\to\infty$
coincides with $((X^j)_{j\in\I}, (L^{jk})_{j,k\in\I})$.

%%%%%%%%%%%%%%%%%%%%%%%%%%%%%%%%%%%%%%
\section{Results}
%%%%%%%%%%%%%%%%%%%%%%%%%%%%%%%%%%%%%%%
\subsection{Existence of a weak solution}
%%%%%%%%%%%%%%%%%%%%%%%%%%%%%%%%%%%%%%%%
We make the following assumption.
 \vskip 3mm \noindent 

\noindent \As[A3]\quad A probability measure $\mu$ on $\mX$ has the log derivative $\mathsf{d}_{\mu } (x,\eta)\in L_{loc}^1(\R^d\times \mX, \mu^{[1]})$, i.e., for any $f\in C_0^\infty(\R^d)\times\mathcal{P}$,
\begin{align*}
&-\int_{\R^d\times \mX}\nabla_xf(x,\eta)
\mu^{[1]}(dx d\eta)
= \int_{\R^d\times \mX}\mathsf{d}_{\mu } (x,\eta)f(x,\eta)\mu^{[1]}(dx d\eta)
\\
&\qquad+\int_{\{(x,\eta): \eta\in \mX, x\in S_\eta\}}\bm{n}_{\eta}(x)f(x,\eta)\mathcal{S}_{\eta}(dx)\rho(x)\mu_{x}(d\eta),
\end{align*}
where $\mathcal{S}_\eta$ is the surface measure on $S_\eta$,
$$
S_\eta=\{x\in\R^d : |x-y|= r \mbox{ for some $y\in\eta$ }\},
$$
and
$\bm{n}_{\eta}(x)$ is the inward normal vector of $S_\eta$ at $x$.

We can extend the notion of the log derivative in distribution and write
$$
\overline{\mathsf{d}}_{\mu } (x,\eta)=\mathbf{1}_{S_\eta}(x)\mathsf{d}_{\mu } (x,\eta)
+\mathbf{1}_{\partial S_\eta}(x)\bm{n}_{\eta}(x)\delta_x.
$$
If the log derivative exists, we put
$b (x,\eta)  = \frac{1}{2} \mathsf{d}_{\mu } (x,\eta)$.
The following result is a modification of \cite[Theorem 26]{o.isde}.

\vskip 3mm

%%%%%%%%%%%%%%%%%%%%
\begin{prp}\label{prp:31}%%%%%%%%%%%%%%%%%%%%%%%%%%%%%%%%%%%%%%
Assume the conditions \As[A0]--\As[A3]. 
There then exists $\mH\subset \mX$ with $\mu(\mH)=1$ such that $(\X= \ml(\Xi), \P_{\xi}), \xi\in\mH$ satisfies the infinite-dimensional stochastic differential equation of Skorohod type (ISKE)
\begin{align}
&dX_t^j = dB_t^j
+ b\bigg(X_t^j, \{X_t^k\}_{k\not=j}\bigg)dt
+ \sum_{k\not=j}\bigg(X_t^j - X_t^k \bigg) dL_t^{jk},
\tag{ISKE}\label{ISKE}
\\ \nonumber
&|X_t^j - X_t^k|\ge r, \quad j,k\in\N, \ j\not=k, \ t\ge 0,
\end{align}
where $L_t^{jk}$, $k,j \in\N$, is a non-decreasing function satisfying
\begin{align}\nonumber
L_t^{jk}=L_t^{kj}=\int_0^t \mathbf{1}(|X_s^j-X_s^k|=r)dL_s^{jk}.
\end{align}
\end{prp}%%%%%%%%%%%%%%%%%%%%%%%%%%%%%%%%%%%%%%%%%%%%%%%%%%%%
%%%%%%%%%%%%%%%%%%%
\noindent{\it Proof.}  Let $m\in \N$.
The coordinate function $x^j$ is locally in the domain of $\mathcal{D}^{\mu^{[m]}}$, and
$X^j$ is thus a Dirichlet process of $(\X^m, \Xi^m)$.
Applying Fukushima decomposition (\cite[Theorem 5.5.1]{fot} 
to $x^j $ yields
\begin{align}\label{:3f}
X_t^{m,j } - X_0^{m,j} = M_t^{[x^j]} + N_t^{[x^j]},\quad \text{ under $\P_{(\x^m,\eta)}^{[m]}$}
.\end{align}
Here, $M^{[x^j]}$ is a martingale additive functional locally of finite energy
and $N^{[x^j]}$ is a continuous additive functional locally of zero energy.
Through a straightforward calculation using \As[A2], we have
for $f \in C_0^\infty (\R^d) \otimes \mathcal{D}_\circ$ that

\begin{align*}
&-\mathcal{E}^{\mu^{[m]}} ( x^i, f)
= \int_{\R^d\times \mX}\mathsf{d}_{\mu } (x^j, \{x^k \}_{k \not=j}^m \cup\eta)f(x,\eta)\mu^{[m]}(dx d\eta)
\\
&\qquad+\int_{\{(x,\eta): \eta\in \mX, x\in S_\eta\}}f(x,\eta)\bm{n}_{\eta}(x^j)
\mathcal{S}_{\eta}(dx^j)\mu_{\x^m}(d\eta)d\x^{m \diamond j},
\end{align*}
where $\x^{m\diamond j} = (x^k)_{k \not= j}^m$.
Hence, by \cite[Theorem 5.2.4]{fot}, we deduce that
\begin{align}\nonumber
N_t^{[x^j]} &= \int_0^t b(X_s^{m,j}, \{X_s^{m,k} \}_{k \not= j}^{m} \cup \Xi_s^{[m]}) ds
+ \sum_{\substack{1 \le k \le m \\ k \not=j} }\int_0^t (X_s^{m,j}- X_s^{m,k} ) dL_s^{m,jk}
\\\label{:3g}
&+ \sum_{k=m+1}^{\infty} \int_0^t (X_s^{m,j} - X_s^{k}) dL_s^{m,jk} 
,\end{align}
where 
$L_t^{m,jk}$, $1\le j \le m$, $k\in \N$ are increasing functions satisfying
\begin{align*}
&L_t^{m,jk}=\int_0^t \mathbf{1}(|X_s^{m,j}-X_s^{m,k}|=r)dL_s^{m,jk},
\quad j,k =1,2,\dots, m,
\\
&L_t^{m,jk}=\int_0^t \mathbf{1}(|X_s^{m,j}-X_s^k|=r)dL_s^{m,jk},
\quad j= 1,2, \dots,m, \ k = m+1, \dots.
\end{align*}
Here, we used the relation $\Xi^m = \{X^j\}_{j=m+1}^\infty$ from the consistency condition \eqref{:22k}.
We put
$$
\D^m [f,g]=
\frac{1}{2}\sum_{i=1}^m \frac{df}{dx_i}\frac{dg}{dx_i} + \mathbb{D}(f,g),
\quad f,g \in C_0^\infty ((\R^d)^m )\otimes \mathcal{D}_\circ.
$$
For $1\le j, \ell \le m$
\begin{align*}
&2\D^m[x^j f, x^j]-\D^m [(x^j)^2,f] = 2\D^m[x^j,x^j]f = f,
\\
&2\D^m[(x^j\pm x^\ell) f, (x^j\pm x^\ell)]-\D^m [(x^j\pm x^\ell)^2,f]
\\
&\qquad\qquad= 2\D^m [x^j\pm x^\ell, x^j\pm x^\ell]f =\begin{cases}
0, \quad &(j=\ell),
\\
2f &(j\not=\ell).
\end{cases}
\end{align*}
Then, from \cite[Theorem5.2.3]{fot},
\begin{align}\label{:3m}
\langle  M^{[x^j]}, M^{[x^\ell]}\rangle_t =\begin{cases}
0, \quad &(j\not=\ell),
\\
t &(j=\ell).
\end{cases}
\end{align}
Combining \eqref{:3f}, \eqref{:3g} and \eqref{:3m} with the consistency \eqref{:22k},
we obtain the proposition.
\qed
%%%%%%%%%%%%%%%%%%%%%%%%%%%%%%%%%%%%%%%%%%%%%%

\vskip 3mm

Let $\mH$ and $\mX_{\rm sde}$ be Borel subsets of $\mX$ such that 
\begin{align}\notag%\label{:A1a}
\mu (\mH) = \mu (\mX_{\rm sde})=1, \quad \mH \subset \mX_{\rm sde} \subset \mX_\infty.
\end{align}
Let $b$ be a measurable function on $\R^d \times \mX$ that has a finite value on
\begin{align*}
\mX_{\rm sde}^{[1]} = \{ x,\eta) \in \R^d \times \mX : \{x\} \cup \eta \in \mX_{\rm sde}\}
.\end{align*}
Let $\ml$ be the label defined by \eqref{:22d}.
We put $\H = \ml (\mH)$ and $\S_{\rm sde}=\ml (\mX_{\rm sde})$. 
Here, $\H$ is the set of initial starting points of solutions and $\S_{\rm sde}$ is the set in which the coefficient of \eqref{ISKE} is well defined.
We consider the following ISKE with \eqref{:A1c} and \eqref{:A1d}:
\begin{align}\tag{ISKE}
&dX_t^j = dB_t^j
+ b\bigg( X_t^j, \{ X_t^k \}_{k\not= j} \bigg)dt
+ \sum_{k\not=j}\bigg(X_t^j - X_t^k \bigg) dL_t^{jk},\quad j\in\N,
\\ \label{:A1c}
&\X_0 \in \H, \quad \X \in W(\S_{\rm sde}),
\\ \label{:A1d}
&L_t^{jk}=\int_0^t \mathbf{1}(|X_s^j-X_s^k|=r)dL_s^{jk}
\quad j, k \in \N,
\end{align}
where $L_t^{jk}, j,k \in \N$ is a non-decreasing real-valued function starting from zero.

\begin{dfn}[weak solution]\label{d:71a}
By a weak solution of \eqref{ISKE} with \eqref{:A1c} and \eqref{:A1d},
we mean an $(\R^d)^N \times \R_+^{\N \times \N} \times (\R^d)^\N$-valued stochastic process $(\X, \L, \B)$ defined on a probability space $(\Omega, \cF, P)$ with a reference family $\{ \cF_t \}_{t\ge 0}$ such that

\noindent {\rm (i)} $(\X, \L)$ is an $\{ \cF_t \}_{t\ge 0}$-adapted $\S_{\rm sde}\times \R_+^{\N\times\N}$-valued process satisfying \eqref{:A1c} and \eqref{:A1d};

\noindent {\rm (ii)} $\B=(B^j)_{j\in\N}$ is an $\R^\N$-valued $\{ \cF_t \}_{t\ge 0}$-Brownian motion with $\B_0={\bf 0}$;

\noindent {\rm (iii)} $\{ b(X_t^j, \{ X_t^k \}_{k\not= j})\}_{j\in\N}$ is a family of $\{ \cF_t \}_{t\ge 0}$-adapted processes with 
$$
E \Big[ \int_0^T |b(X_t^j, \{ X_t^k \}_{k\not= j})|dt
\Big]< \infty \quad \text{for all $T$}; and
$$

\noindent {\rm (iv)} with probability one, $(\X, \L, \B)$ satisfies for all $t \ge 0$ that
$$
X_t^j = X_0^j + B_t^j
+ \int_0^tb\bigg( X_u^j, \{ X_u^k \}_{k\not= j} \bigg)du
+ \sum_{k\not=j}\int_0^t\big(X_u^j - X_u^k \big) dL_u^{jk},\quad j\in\N.
$$
\end{dfn}
%%%%%%%%%%%%%%%%%%%%%%%%%%%%%%%%%%%%%%%%%%%%%%%%%

\begin{rem}\label{R:1}%%%%%%%%%%%%%%%%%%%%%%%%%%%%%%%%%%%%%%%%%%%%%%%%%%
Let $\mu$ be a canonical Gibbs state with the potentials $(\Phi, \Psi = \Psi_{\rm hard}+\Psi_{\rm sm})$ such that $\Phi$ and  $\Psi_{sm}$ are smooth. 
From the same argument as \cite[Lemma 13.5]{o-t.tail} 
the log derivative of $\mu$ exists and is represented as
\begin{align}\nonumber%\label{:31w}
\mathsf{d}_{\mu } (x,\eta) = -\nabla \Phi(x) -\sum_{y\in \eta}\nabla \Psi_{\rm sm}(x-y).
\end{align}
It follows from applying Proposition \ref{prp:31} that 
$(\X_t, \P_{\x})$ is a weak solution (SDE-$\N$).
\end{rem}%%%%%%%%%%%%%%%%%%%%%%%%%%%%%%%%%%%%%%%%%%%%%%%

%%%%%%%%%%%%%%%%%%%%%%%%%%%%%%%%%%%%%%%%%%%%%%%%%%
\subsection{ Statement of the results}
%%%%%%%%%%%%%%%%%%%%%%%%%%%%%%%%%%%%%%%%%%%%%%%%%%

We study the existence and uniqueness of strong solutions of \eqref{ISKE} with \eqref{:A1c} and \eqref{:A1d},
whose definition are given in Definitions \ref{d:74a} and \ref{d:76a}. 
In \cite{o-t.tail}, we developed a general theory of the existence of a strong solution and the pathwise uniqueness of solutions for ISDEs concerning interacting Brownian motions. 
We apply the argument made in the cited paper.

Let $(\X, \B)$ be an $(\R^d)^\N \times (\R^d)^\N$-valued continuous process defined on a standard filtered space $(\Omega, \mathcal{F}, P, \{\mathcal{F}_t\})$. 
The regular conditional probability $P_\x = P(\cdot | \X_0=\x)$ then exists for $P\circ \X_0^{-1}$-a.s. $\x$.

We cite the conditions for a weak solution of $(\X, \B, P)$ given in \cite[Section 3.2]{o-t.tail}.

\vskip 3mm
\noindent \As[SIN] \quad $P( \X \in W(\S_{hard}))=1$.

\vskip 3mm
\noindent \As[$\mu$-AC] ($\mu$-absolutely continuity condition)
\begin{align}\nonumber
P(\mfu(\X_t)\in \cdot) \prec \mu \quad \ \text{for all $t>0$}
,\end{align}
where for the two Radon measures $m_1$ and $m_2$, 
we write $m_1 \prec m_2$ if $m_1$ is absolutely continuous with respect to $m_2$.
\vskip 3mm
\noindent \As[NBJ] (No big jump condition) \quad 
 $\forall \ell , \forall T \in \mathbb{N} $
$$ 
P ( \mathsf{m}_{r,T}(\X) < \infty  ) = 1, 
$$
where
 \begin{align}\nonumber
 &\mathsf{m}_{\ell,T}(\X)
  = \inf \{ m \in \mathbb{N}\, ;\, 
|X ^n (t)|>\ell,  \forall n  > m, \forall t\in [0,T] \} 
.\end{align}

For a topological space $S$, we denote by $\mathcal{B}(S)$ the topological Borel field of $S$. 
We say a family of strong solutions $\X=F_\x (\B)$ starting at $\x$ for $P\circ \X_0^{-1}$-a.s. $\x$ satisfies the measurable family condition if
\vskip 3mm
\noindent \As[MF] 
 \quad $P(F_\x(\B) \in A)$ is $\overline{\cB((\R^d)^\N)}^{\X_0}$-measurable for any $A \in \cB(W(\R^d))$,
\vskip 3mm
\noindent where $\overline{\cB((\R^d)^\N)}^{\X_0}$ is the completion of $\B(W((\R^d)^\N ))$ with respect to $P\circ \X_0^{-1}$.

The tail $\sigma$-field on $\mX$ is defined as
$$
\mathcal{T}(\mX) = \bigcap_{r=1}^{\infty} \sigma(\pi_{S_r^c}).
$$
We introduce the following condition on a probability measure $\mu$.
\vskip 3mm \noindent 
\As[TT] (tail trivial) \quad
$\mu (\mathsf{A}) \in \{ 0, 1\}$ \quad for any $\mathsf{A} \in \mathcal{T} (\mX)$.

\vskip 3mm

We make the following assumption.

\vskip 3mm 
\noindent 
\As[R] $\Phi \equiv 0$ and $\Psi_{\rm sm}$ is a Ruelle's class potential such that 
\begin{align}\notag %\label{:33f}
\sup_{\xi\in\mX}\sum_{x\in\xi}|\nabla\Psi_{\rm sm}(x)|<\infty, \quad
\sup_{\xi\in\mX}\sum_{x\in\xi} |\nabla^2 \Psi_{\rm sm}(x)| <\infty.
\end{align}

For a subset $\I$ of $\N$, we put $\I^c = \N\setminus \I$.
We introduce another condition \As[$\mathcal{I}$-IFC] weaker than \As[IFC] given in Section \ref{s:52}.
Let $T\in \N$ and $\M =\{(\I_i,t_i)\}_{i=0}^M$ be sequences of pairs of finite index sets and times such that
\begin{align}\notag %\label{:32d}
\I_0 \supset \I_1 \supset \cdots \supset \I_{M-1} =: \I_*,
\quad
t_i = \frac{iT}{M}, \quad i=0,1, \dots,M-1.
\end{align}
We introduce the sequence of SDEs 
\begin{align*} \nonumber %%%%%%%%%%%%%%%%
&dY_t^{\I_i,j} =  dB_t^j + 
b^{\I_i,j}_{\X} (Y_t^{\I_i,j},\mathbf{Y}_t^{\I_i})dt
+\sum_{k\in \I_i \setminus \{j\}}(Y_t^{\I_i,j} - Y_t^{\I_i,k}) dL_t^{\I_i,jk}
 \\ \tag{SKE$_{\X}$($\M$)} \label{SKXIi}
&\qquad+\sum_{k\in \I_i^c}^\infty (Y_t^{\I_i,j} - X_t^{k}) dL_t^{\I_i,jk},
\quad j\in \I_i, \quad t \in [t_i, t_{i+1}],
\\
&Y_{t_0}^{\I_0,j} = x^j, j\in \I_0, \quad Y_{t_i}^{\I_i,j} = Y_{t_i}^{\I_{i-1}},
\ j\in \I_i, \ i=1,2,\dots, M-1,
\end{align*}
where $i= 0,1,\dots,M-1$, $\displaystyle{b_{\X}^{\I_i,j}(t,\y)=b\bigg(y^j,\{y^k\}_{k\in\I_i \setminus \{j\}}+ \{X_t^k\}_{k\in \I_i^c}\bigg)}$, and $L_t^{\I_i,jk}$, $j\in\I_i$, $k\in \N$ are increasing functions satisfying
\begin{align*}
&L_t^{\I_i,jk}=\int_{t_i}^{t_{i+1}} \mathbf{1}(|Y_s^{\I_i,j}-Y_s^{\I_i,k}|=r)dL_s^{\I_i, jk},
\quad j,k \in \I_i,
\\
&L_t^{\I_i,jk}=\int_{t_i}^{t_{i+1}} \mathbf{1}(|Y_s^{\I_i,j}-X_s^k|=r)dL_s^{\I_i, jk},
\quad j\in \I_i, \ k\in \I_i^c.
\end{align*}
.

We introduce the sequence $\{\Lambda(\I_i, [t_i,t_{i+1}])\}_{i=0}^{M-1}$ of the events defined by 
$$
\Lambda(\I_i, [t_i, t_{i+1}]) = \{ |X_u^{k} - Y_u^{\I_i, j}|>r, \quad u \in [t_i, t_{i+1}], \ j\in \I_i, \ k \in \I_i^c \}
$$
and put $\Lambda^\M = \bigcap_{i=0}^{M-1} \Lambda(\I_i, [t_i, t_{i+1}])$.

We denote by $\mathcal{C}^{\I_0, \I_*^c}$ the completion of $\cB(W_{\bf 0}((\R^{d})^{\I_0})\times W((\R^d)^{\I_*^c}))$
with respect to $P_{\x}\circ (\B^{\I_0}, \X^{\I_*^c})^{-1}$.
Let $\u\in W ((\R^d)^\N)$ and $(\mathbf{v}, \w)\in W_{\bf 0}((\R^{d})^{\I_0})\times W((\R^d)^{\I_*^c})$.
We put $\cB_t(W((\R^d)^\N))=\sigma[\u_s : 0\le s \le t]$ and denote
by $\mathcal{C}_t^{\I_0, \I_*^c}$ the completion of $\sigma[(\mathbf{v}_s, \w_s) : 0\le s \le t]$
with respect to $P_{\x}\circ (\B^{\I_0}, \X^{\I_*^c})^{-1}$.
We then make the following assumption.
\vskip 3mm
\noindent \As[$\mI$-IFC] \quad For each $\M= \{ (\I_i, t_i) \}_{i=0}^M$, 
the pathwise uniqueness of solutions of (\ref{SKXIi}) on $\Lambda^\M$ holds
under $P_{\x}$ for $P\circ \X_0^{-1}$-a.s $\x$. Moreover, 
there exists a $\mathcal{C}^{\I_0, \I_*}$-measurable function 
$$
F_\x^\M : W_0((\R^d)^{\I_0}) \times W((\R^d)^{\I_*^c} ) \to W((\R^d)^\N)
$$
such that $F_\x^\M$ is $\mathcal{C}_t^{\I_0, \I_*}/ \cB_t(W((\R^d)^\N))$-measurable for each $t$ 
and satisfies
$$
F_\x^\M(\B^{\I_0}, \X^{\I_*^c})_s = (Y_s^{\I_i}, \X_s^{\I_i^c}), \ s\in [t_i, t_{i+1}], \
i=0,\dots,M-1,
$$
on $\Lambda^{\M}$ under $P_{\x}$ for $P\circ \X_0^{-1}$-a.s $\x$.
\vskip 3mm
We also introduce the other condition \As[FCP].
We first introduce measurable subsets of $W(\S_{\rm hard})$. 
For $\varepsilon >0$, $0 \le s < t < \infty$ and a bounded connected open set $\C$ of $\R^d$, 
we denote by $\mathfrak{C}(\varepsilon, [s,t], \C)$ the set of all elements $\X= (X^1, X^2, \dots)$ of $W(\S_{\rm hard})$ such that
\begin{align}\notag
& U_{(r+\varepsilon)/2}(X_u^j) \in \C, \ u \in [s,t], \quad\text{if $X_s^j \in \C$,}
\\ \notag
& U_{(r+\varepsilon)/2}(X_u^j) \in \R^d \setminus \C, \ u \in [s,t], \quad\text{if $X_s^j \in \R^d\setminus \C$}
.\end{align}
%Here $U_\alpha (A)$, $\alpha >0$, denotes the open $\alpha$-neighbourhood of a set $A\subset \R^d$, and $U_\alpha (x)$ is the abbreviated form of $U_\alpha (\{x\})$. 
For $\varepsilon >0$, $p, T, a, M \in \N$, we denote by $\mathfrak{C}(\varepsilon, p, T, a, M)$ the set of elements $\X$ of $W(\S_{\rm hard})$ such that
%\begin{align}
$\displaystyle{
\X \in \bigcap_{i=0}^{M-1}\mathfrak{C} (\varepsilon, [t_i, t_{i+1}], \C_i),
}$
%\end{align}
for $t_i = \frac{iT}{M}$, $i=0,1,\dots,M-1$,
and some decreasing sequence $\bbC=\{\C_i\}_{i=0}^{M-1}$ of open subsets of $\R^d$ with
\begin{align}\nonumber %
\C_0 \subset U_{a+M + M^p}(0), \ U_{\varepsilon}(\C_{i+1}) \subset \C_i, \ 
 0\le i\le M-2 
 \\ \label{;41c}
 \text{ and, } U_{a+M}(0) \subset \C_{M-1}.
\end{align}
We denote a measurable subset $\mathfrak{C}$ of $W (\S_{\rm hard})$ by
\begin{align}\nonumber
\mathfrak{C} = \bigcup_{\varepsilon >0} \bigcup_{p=1}^\infty 
\bigcap_{T=1}^\infty \bigcap_{a=1}^\infty \bigcap_{M_0=1}^\infty \bigcup_{M = M_0}^{\infty}
\mathfrak{C}(\varepsilon, p, T, a, M)
.\end{align}
Note that $\X \in \mathfrak{C}$ implies $\theta_t\X \in \mathfrak{C}$ for any $t>0$,
where $\theta_t \X(u)= \X(u+t)$. 
We remark that we can define $\mathfrak{C}$ by a countable collection of bounded open subsets of $\R^d$,
because if $\X \in \mathfrak{C}(\varepsilon, [s,t], \C)$, there exists $\varepsilon' > 0$ and a polyhedron $\C'$ with vertices in $\varepsilon'\Z^d$ such that $\X \in \mathfrak{C}(\varepsilon', [s,t], \C')$. 
We then assume $\bbC$ is chosen from a countable family $\mathcal{A} = \{ \bbC (\ell) \}_{\ell\in\N}$.

We make the following assumption, called the finite cluster property \As[FCP].

\vskip 3mm

\noindent \As[FCP] \quad $ P ( \X \in \mathfrak{C} )=1$.

\vskip 3mm

The main theorem of this article is the following.

\begin{thm}\label{th:main}%%%%%%%%%%%%%%%%%%%%%%%%%%%%%%%%%
Assume \As[TT]. 

\noindent{\rm (i)} Assume \As[A0], \As[A2], \As[A3], and \As[R]. 
Put $(\X, P) = (\ml_{\rm path}(\Xi), \P_\mu)$. 
Then, for $\mu\circ \ml^{-1}$-a.s. $\x$,
$(\X, P_\x)$ is a strong solution of (\ref{ISKE}) with \eqref{:A1c} and \eqref{:A1d} starting at $\x$.
Moreover, $(\X, P)$ satisfies \As[MF], \As[$\mI$-IFC], \As[FCP], \As[$\mu$-AC], \As[SIN], and \As[NBJ].

\noindent{\rm (ii)}
\eqref{ISKE} with \eqref{:A1c} and \eqref{:A1d} has a family of unique strong solutions $\{F_\x\}$ starting at $\x$ for $P\circ \X_0^{-1}$-a.s. $\x$ under the constraints of 
\As[MF], \As[$\mI$-IFC], \As[FCP], \As[$\mu$-AC], \As[SIN], and \As[NBJ].
\end{thm} %%%%%%%%%%%%%%%%%%%%%%%%%%%%%%%%%

\begin{rem} If \As[A1] is satisfied, we can decompose $\mu$ as 
$$
\mu = \int_{\mX}\mu(d\eta) \mu_{{\rm Tail}}^{\eta},
$$
where $\mu_{{\rm Tail}}^{\eta}= \mu(\cdot | \mathcal{T}(\mX))(\eta)
$ is the regular conditional distribution with respect to the tail $\sigma$-field.
Note that \As[TT] for $\mu_{{\rm Tail}}^{\eta}$ holds.
In the case that $\mu$ satisfies \As[R], \As[MF], \As[$\mI$-IFC], \As[FCP], \As[SIN], and \As[NBJ],
$\mu_{{\rm Tail}}^{\eta}$ also does for $\mu$-a.s. $\eta$.
Hence, assuming \As[$\mu_{\rm Tail}$-AC] for $\mu$-a.s. $\eta$ instead of \As[$\mu$-AC], 
the counterpart of Theorem \ref{th:main} is derived.
The constraint of \As[$\mu_{\rm Tail}$-AC] means that 
there is no $A\in \mathcal{T}(\mM)$ such that
for $\mu_{{\rm Tail}}^{\eta}$-a.s. $\xi$, 
$$
P(\mfu(\X_s)\in A | \mfu(\X_0)=\xi) \not= P(\mfu(\X_t)\in A | \mfu(\X_0)=\xi)
$$
for some $0\le s<t$. It is then possible that another solution $\X'$ that changes the tail exists. 
See \cite[Section 3.3]{o-t.tail}.
\end{rem}
\vskip 3mm

\begin{exa} Recalling Remark \ref{R:1}, we present two examples for Theorem \ref{th:main}.
\begin{itemize}
\item[{\rm (i)}] 
Lennard-Jones 6-12 potential
$( d= 3$, 
$\Psi_{\rm sm}(x) = \Psi_{6,12}(x) = |x|^{-12}-|x|^{-6})$
\begin{align} \nonumber
b(x, \{y^k\})= \frac{\beta }{2} 
\sum _{k} \Big\{
\frac{12(x-y^k)}{|x-y^k|^{14}} - 
\frac{6(x-y^k)}{|x-y^k|^{8}}\, \Big\}.
\end{align}

\noindent \item[\rm (ii)] Riesz potentials
$( d < a \in \N  $ and $\Psi_{\rm sm}(x) = \Psi_a(x)=(\beta / a )|x|^{-a})$
\begin{align}\nonumber
b(x, \{y^k\})= \frac{\beta }{2}
\sum _{k} \frac{x-y^k}{|x-y^k|^{a+2}}.
\end{align}
\end{itemize}
\end{exa}

%%%%%%%%%%%%%%%%%%%%%%
\section{Proof of the main theorem}
%%%%%%%%%%%%%%%%%%%%%%

We prepare lemmas for proving the main theorem.

\subsection{Finite cluster property}

%%%%%%%%%%%%%%%%%%%%%%%%%%%%%%%%%%%%%%
 \begin{lmm}[Lemma 2.6 in \cite{FRT}]\label{l:41b}
 Assume \As[A0], \As[A2], and \As[R]. 
 Then, \As[FCP] for $\X= \ml_{\rm path}(\Xi)$ holds.
 \end{lmm}

\begin{rem}
{\rm (i)} In the proof of Lemma 2.6 in \cite{FRT}, we used an estimate of the modulus of continuity of each ball derived through  Lyons--Zhen decomposition \cite[Section 5.7]{fot}.
We again use the estimate here.
\\
\noindent
{\rm (ii)} 
In the proof, we used a property of the continuum percolation model associated with $\mu$. See \cite[Lemma 2.4]{FRT}.  It is seen that this property holds under \As[R] but not under \As[A1], which is an obstacle to generalizing Theorem \ref{th:main} for quasi-Gibbs states. 
It is an interesting problem to study percolation theory for quasi-Gibbs states. See, for instance, Ghosh \cite{G16}.
 
%{\rm (ii)}  In the proof of Lemma 2.6 in \cite{FRT} we used the estimate of the cluster size $\sharp \C$ under the distribution $\mu = \mu^{\Psi}$ in \cite[Lemma 2.5]{FRT}. We used the stretched exponential decay of $\Psi_{\rm sm}$ for the estimate. In this article, we can omit the condition.
\end{rem}

\subsection{ On the Lipschitz continuity of $b_{\X}^{\I}$}

In this subsection, we examine the Lipschitz continuity of $b_{\X}^{\I}$ for a finite subset $\I$ of $\N$.
We assume \As[A0] and \As[R]. We recall that
$b(x,\eta) = -\frac{1}{2}\sum_{y\in\eta} \nabla \Psi_{\rm sm}(x-y)$ in Remark \ref{R:1}. 
Let $\I$ be a finite subset of $\N$. 
We introduce the domain of the configurations of balls indexed by $\I$ given by
$$
D^\I= \{\x^\I\in (\R^d)^\I:
|x^j-x^k|>r, \quad j\not=k, \ j,k \in \I\}.
$$
For $\x=(x^k)_{k\in\I} \in D^\I$ and $\eta \in \mX$ with 
$\{ x^j \}_{j\in\I} \cap \eta = \emptyset, \ \{ x^j \}_{j\in\I} \cup \eta \in \mX$,
we set
$$
b^\I (\x, \eta) = (b(x^j, \{x^k\}_{k\in \I \setminus \{j\}}+\eta))_{j\in\I}.
$$
Let $K(\eta)\in [0,\infty]$ be a function defined by
\begin{align}\nonumber
K(\eta)&= \sup \Big\{ \frac{b^\I(\x,\eta)-b^\I(\y,\eta)}{|\x-\y|} :\x\not= \y, 
\x, \y \in D^\I
\\  \nonumber
&\qquad \{ x^j \}_{j\in\I} \cap \eta = \{ y^j\}_{j\in\I} \cap \eta = \emptyset,
\ \{ x^j \}_{j\in\I} \cup \eta, \{ y^j\}_{j\in\I} \cup \eta \in \mX
\Big\}.
\end{align}
It follows from \As[R] that 
\begin{align}\label{:42e}
K:= \sup_{\eta\in\mX}K(\eta) < \infty
\end{align}
and it follows from
$$
b_{\X}^{\I}(t, \x) = b^{\I}\big(\x, \mfu (\X_t^{\I^c})\big)
$$
that
\begin{align}\label{:42g}
|b_{\X}^{\I}(t, \x) - b_{\X}^{\I}(t, \y)| \le K |\x - \y|.
\end{align}

\begin{rem}
In \cite[Section 11]{o-t.tail}, the Lipschitz continuity of $b_{\X}^{\I}$ was discussed in the case that $\mu$ is a quasi-Gibbs state. 
The method presented in that section is applicable to the case with a hard core.  
\end{rem}

%%%%%%%%%%%%%%%%%%%%%%%%%%%%%%%
\subsection{$\mI$-IFC}

\begin{lmm}\label{l:43a}
Assume \As[A0], \As[A2], \As[A3], and \As[R]. Let $(\X, P_\x) = (\ml_{\rm path}(\Xi), \P_{\mfu(\x)})$.
\As[$\mI$-IFC] then holds for $\X$.
\end{lmm}

\noindent {\it Proof.} 
Let $\M= \{(\I_i,t_i)\}_{i=0}^{M-1}$.
Suppose that $\omega \in \displaystyle{\Lambda^\M}$.
Then, $\{ \Y^{\I_i}\}$ satisfies
\begin{align} \nonumber %%%%%%%%%%%%%%%%
&dY_t^{\I_i,j} =  dB_t^j + 
b^{\I_i,j}_{\X} (t,\mathbf{Y}_t^{\I_i})dt
+\sum_{k\in \I_i \setminus \{j\}}(Y_t^{\I_i,j} - Y_t^{\I_i,k}) dL_t^{\I_i,jk},
 \\ \label{SKXI_Lambda}
&\qquad\qquad\qquad\qquad\qquad
\quad j\in \I_i, \quad t \in [t_i, t_{i+1}], \quad i=0, 1, \dots, M-1.
\end{align}
Let $\Y^{\I_0}$ and $\tilde{\Y}^{\I_0}$ be solutions of \eqref{SKXI_Lambda} with $i=0$.
Put
\begin{align*}
&w_t = \B_t^\I + \int_0^t b^{\I_0,j}_{\X} (s,\mathbf{Y}_s^{\I_0})ds, \quad t\in [0,t_1],
\\
&\tilde{w}_t = \B_t^\I + \int_0^t b^{\I_0,j}_{\X} (s,\tilde{\mathbf{Y}}_s^{\I_0})ds, \quad t\in [0,t_1].
\end{align*}
It then follows from the Lipschitz continuity \eqref{:42g} and \eqref{:42e} that
\begin{align*}
\|w - \tilde{w}\|_t \le \int_0^t |b^{\I_0,j}_{\X} (s,\mathbf{Y}_s^{\I_0}) 
- b^{\I_0,j}_{\X} (s,\tilde{\mathbf{Y}}_s^{\I_0})|ds
\le K \int_0^t |\mathbf{Y}_s^{\I_0}-\tilde{\mathbf{Y}}_s^{\I_0}|ds.
\end{align*}
$\triangle_{0,t_1,\delta}(w)<\infty, \ \triangle_{0,t_1,\delta}(\tilde{w})<\infty,
\ \sup_{0\le s\le t_1}|w|<\infty, \ \sup_{0\le s\le t_1}|\tilde{w}|<\infty$,
and we thus have the pathwise uniqueness in the case that $i=0$ from \eqref{:23c} and \eqref{:23d}. 
Repeating this procedure, we obtain the pathwise uniqueness for $i=1,2,\dots, M-1$. 
Using an approximation process as in \eqref{SKE(n)_I} and the above estimate,
we can show the existence of a strong solution $F_\x^{\M}$.
\qed

%%%%%%%%%%%%%%%%%%%%%%%%%%%%%%%
\subsection{Proof of Theorem \ref{th:main} }

We use the following lemma, which is a modification of \cite[Theorem 3.1]{o-t.tail}, 
with the condition \As[IFC] replaced by the pair of conditions \As[$\mI$-IFC] and \As[FCP].

\begin{lmm}\label{l:44a}
Assume \As[TT] for $\mu$. Assume that \eqref{ISKE} with \eqref{:A1c} and \eqref{:A1d} has a weak solution $(\X, \B)$ under $P$ satisfying \As[$\mI$-IFC], \As[FCP], \As[$\mu$-AC], \As[SIN], and \As[NBJ]. 
Then, \eqref{ISKE} with \eqref{:A1c} and \eqref{:A1d} has a family of unique strong solutions $\{F_\x\}$ starting at $\x$ for $P\circ \X_0^{-1}$-a.s. $\x$ under the constraints of 
\As[MF], \As[$\mI$-IFC], \As[FCP], \As[AC] for $\mu$, \As[SIN], and \As[NBJ].
\end{lmm}
From this lemma, Theorem\ref{th:main} is shown, 
if we check that the weak solution $(\ml_{\rm path}(\Xi), \P_{\mu})$ of \eqref{ISKE} 
satisfies $\As[MF]$, \As[$\mI$-IFC], \As[FCP], \As[$\mu$-AC], \As[SIN], and \As[NBJ]. 
\As[MF] is obvious. 
\As[AC] is derived from the reversibility of the process $\Xi$ with respect to $\mu$.
\As[SIN] and \As[NBJ] follows from $\Psi_{\rm sm}$ being of Ruelle's class with a hard core.
See Lemmas 10.2 and 10.3 in \cite{o-t.tail}. 
\As[$\mI$-IFC] is derived from Lemma \ref{l:43a}.
Hence, it is enough to show Lemma \ref{l:44a} to prove Theorem\ref{th:main}.

In the proof \cite[Theorem 3.1]{o-t.tail}, \As[IFC] is used in \cite[Lemma 4.2]{o-t.tail}.
That is to say the existence of a function
$F_{\x}^\infty : W_0((\R^d)^\N)\times  W((\R^d)^\N) \to W((\R^d)^\N)$
satisfying
\begin{itemize}
\item[(i)]  $F_\x^\infty (\B, \X) =\X$ $P_\x$-a.s. and
\item[(ii)] 
$F_\x^\infty (\b, \cdot)$ is $\mathcal{T}_{\rm path}((\R^d)^\N)_{\x, \b}$-measurable for $P_{Br}^\infty := P\circ B^{-1}$- a.s $\b$,
\end{itemize}
where $\mathcal{T}_{\rm path}((\R^d)^\N)_{\x, \b}$ is the completion of
$$
\mathcal{T}_{\rm path}((\R^d)^\N) := \bigcap_{\I \subset \N, \sharp\I < \infty} \sigma (\X^{\I^c})
$$
with respect to $P_\x \circ (\X,\B)^{-1}( \cdot | \B=\b)$.
We  construct $F_\x^\infty$ under the conditions \As[SIN], \As[$\mI$-IFC], and \As[FCP].
Let $\bbC =\{\C_i\}_{i=0}^{M-1}$ with \eqref{;41c}. Put
$$
\I(\C_i) =\{j\in \N: X^j_{t_i} \in \C_i\}, \quad
\M (\bbC) = \{(\I(\C_i), t_i) \}_{i=0}^{M-1}
$$
and
$$
F_{\x}^{[\bbC]} = F_{\x}^{\M(\bbC)}
:= \sum_{\M} F_{\x}^{\M}\mathbf{1}(\M(\bbC )=\M)\mathbf{1}_{\Lambda^\M}
 \ \text{ on \ $\bigcup_{\M}\Lambda^\M$}.
$$
Let $a\in \N$. 
From \As[FCP] and \As[$\mI$-IFC], for $P_{Br}^{\infty}$-a.s. $\b$ and $P(\cdot|B=b)$ -a.s $\X$, 
there exists $\bbC \in \mathcal{A}:=\{ \bbC (\ell)\}$ such that $U_{a+M} \subset \C_{M-1}$ and 
$$
F_\x^{[\bbC]}(\b, \X) \in \Lambda^{\M(\bbC)}.
$$
We put 
$\ell (\b, \X) = \min\{ \ell \in \N : F_\x^{\bbC(\ell)}(\b, \X) \in \Lambda^{\M(\bbC(\ell))} \}$
and
$$
F_\x^{(a)}(\b, \X)=\sum_{\ell\in N} \mathbf{1}(\ell=\ell(\b,\X)) F_\x^{[\bbC (\ell)]}(\b, \X).
$$
Let $\I$ be a finite subset of $\N$.
Then,
$F_\x^{(a)}(\b, \cdot)\mathbf{1}(\I(\C_{M-1}) \supset \I)$ is $\overline{\sigma(\X^{\I^c})}$-measurable, where 
$\overline{\sigma(\X^{\I^c})}$ is the completion of $\sigma(\X^{\I^c})$ with respect to $P_\x \circ (\X,\B)^{-1}( \cdot | \B=\b)$.
From \As[SIN], we see that
$\displaystyle{\lim_{a\to\infty}} P_\x (\I(\C_{M-1}) \supset \I | \B=\b)=1$. 
Putting 
$$
F_\x^\infty = \lim_{a\to\infty}F_\x^{(a)},
$$
we see that $F_\x^\infty (\b, \cdot)$ satisfies (ii). 
The claim (i) is derived from \As[$\mI$-IFC]. 
Therefore, Lemma \ref{l:44a} is proved adopting the same procedure used in \cite[Theorem 3.1]{o-t.tail}.
\qed

\vskip 3mm
%\newpage

%%%%%%%%%%%%%%%%%%%%%%%%%%%%%
\section{Appendix}
%%%%%%%%%%%%%%%%%%%%%%%%%%%%

\subsection{Solutions of \eqref{ISKE}}

We give precise definitions of solutions of
the \eqref{ISKE}.

%\begin{dfn}[uniqueness in law]. \label{d:72a}
%We say that the uniqueness in weak solution of \eqref{ISKE} with \eqref{:A1c} and \eqref{:A1d} holds
%if whenever $\X$ and $X'$ are two solutions whose initial distributions coincide,
%then the law of $X$ and $X'$ do as well.
%\end{dfn}

%\begin{dfn}[pathwise uniqueness].\label{d:73a}
%We say that the pathwise uniqueness of weak solution of \eqref{ISKE} with \eqref{:A1c} and \eqref{:A1d} holds
%if, whenever $\X$ and $X'$ are two solutions defined o the same probability space $(\Omega, \cF, P)$ with same reference family $\{\cF_t\}_{t\ge 0}$ and same $(\R^d)^{\N}$-valued Brownian motion $\{\cF_t\}_{t\ge 0}$-Brownian motion $\B$ such that $\x_0 = \X'_0 \in \H$ a.s., then
%$X_t = X'_t$, $t\ge 0$ a.s.
%\end{dfn}
%\vskip 3mm

%Let $\cB (W((\R^d)^\N)$ be the Borel $\sigma$-field on $W((\R^d)^{\N})$and $\cB_t (W((\R^d)^\N) = \sigma [\w_s ; 0\le s \le t]$.
%We set $\cB_t (W_0((\R^d)^\N)$ similarly. 
Let $\overline{\cB}$ and $\overline{\cB}_t$ be the completions of 
$\cB (W((\R^d)^\N)$ and $\cB_t (W((\R^d)^\N)$ with respect to 
$P_{Br}^{\infty}= P(\B \in \cdot)$, respectively.

\vskip 3mm

\begin{dfn}[strong solutions starting at $\x$]\label{d:74a}
A weak solution $\X$ of \eqref{ISKE} with \eqref{:A1c} and \eqref{:A1d} and an $(\R^d)^{\N}$-valued $\{\cF_t\}_{t\ge 0}$-Brownian motion $\B$ on $(\Omega, \cF, P, \{ \cF \}_{t\ge 0})$ is called a strong solution starting at $\x$ if $\X_0=\x$ a.s. and if there exists a function $F_\x : W_0((\R^d)^\N) \to W((\R^d)^\N)$ 
such that $F_\x$ is $\overline{\cB}/ \cB(W((\R^d)^\N)$-measurable, 
and $\overline{\cB}_t/ \cB_t(W((\R^d)^\N)$-measurable for each $t$
and that $F_\x$ satisfies
$\X = F_\x (\B)$ a.s.
We also call $\X= F_\x (\B)$ a strong solution starting at $\x$.
Additionally, we call $F_\x$ itself a strong solution starting at $\x$.
\end{dfn}

\vskip 3mm

\begin{dfn}[unique strong solution starting at $\x$]\label{d:75a}
We say \eqref{ISKE} with \eqref{:A1c} and \eqref{:A1d} has a unique strong solution starting at $\x$ 
if there exists a function $F_\x : W_0((\R^d)^\N) \to W((\R^d)^\N)$ such that,
for any weak solution $(\hat{\X}, \hat{\L}, \hat{\B})$ of \eqref{ISKE} with \eqref{:A1c} and \eqref{:A1d}
$\hat{\X} = F_\x (\hat{\B})$  a.s.
and if, for any $(\R^d)^\N$-valued 
$\{\cF_t\}_{t\ge 0}$-Brownian motion $\B$ defined on $(\Omega, \cF, P, \{\cF\}_{t\ge 0})$ 
with $\B_0=0$, the process $\X=F_\x (\B)$ is a strong solution of \eqref{ISKE} with \eqref{:A1c} and \eqref{:A1d} starting at $\x$.
We also call $F_\x$ a unique strong solution starting at $\x$.
\end{dfn}

We next present a variant of the notion of a unique strong solution.

\begin{dfn}[a unique strong solution under a constraint]\label{d:76a}
For a condition \As[Cond], we say \eqref{ISKE} with \eqref{:A1c} and \eqref{:A1d} has a unique strong solution starting at $\x$ under the constraint of \As[Cond] 
if there exists a function $F_\x : W_0((\R^d)^\N) \to W((\R^d)^\N)$ such that for any weak solution $(\hat{\X}, \hat{\B})$ of \eqref{ISKE} with \eqref{:A1c} and \eqref{:A1d} starting at $\x$ satisfying \As[Cond], it holds that
$\hat{\X} = F_\x (\hat{\B})$ a.s.
and if, for any $(\R^d)^\N$-valued 
$\{\cF_t\}_{t\ge 0}$-Brownian motion $\B$ on $(\Omega, \cF, P, \{\cF\}_{t\ge 0})$ 
with $\B_0=0$, the process $\X=F_\x (\B)$ is a strong solution of \eqref{ISKE} with \eqref{:A1c} and \eqref{:A1d} starting at $\x$ satisfying \As[Cond]. 
We also call $F_\x$ a strong solution starting at $\x$ under the constraint of \As[Cond].
\end{dfn}

For a family of strong solutions $\{F_\x\}$ satisfying \As[MF], we put
$$
P_{\{F_\x\}} = \int P(F_\x(\B)\in \cdot) P \circ \X_0^{-1}(d\x).
$$
Let $(\X, \L, \B)$ be a solution of \eqref{ISKE} with \eqref{:A1c} and \eqref{:A1d} under P.
Suppose that $(\X, \B)$ is a unique strong solution under $P_\x$ for $P\circ \X_0^{-1}$-a.s.
$\x$. 
Let $\{F_\x\}$ be a family of the unique strong solution given by $(\X, \B)$ under $P_\x$. 
Then, \As[MF] is automatically satisfied and $P_{\{F_\x \}}= P\circ \X^{-1}$.

\vskip 3mm

\begin{dfn}[a family of unique strong solutions under constraints] \label{d:77a}
For a condition \As[Cond], we say \eqref{ISKE} with \eqref{:A1c} and \eqref{:A1d} has a family of unique strong solutions $\{F_\x\}$ starting at $\x$ for $P\circ \X_0^{-1}$-a.s. $\x$
under the constraints of \As[MF] and \As[Cond] 
if $\{F_\x\}$ satisfies \As[MF] and $P_{\{F_\x\}}$ satisfies \As[Cond].
Furthermore, (i) and (ii) are satisfied.

\noindent {\rm (i)} For any weak solution $(\hat{\X}, \hat{\B})$ under $\hat{P}$ of \eqref{ISKE} with \eqref{:A1c} and \eqref{:A1d} with $\hat{P}\circ \X_0^{-1} \prec P\circ \X_0^{-1}$ satisfying \As[Cond], it holds that, for $\hat{P}\circ \X_0^{-1}$-a.s. $\x$, $\hat{\X}=F_\x(\hat{\B})$ $\hat{P}_\x$-a.s.,
where $\hat{P}_\x = \hat{P}(\cdot | \hat{\X}_0 =\x)$.

\noindent {\rm (ii)} For an arbitrary $(\R^d)^\N$-valued $\{\cF_t\}$-Brownian motion $\B$ on $(\Omega, \cF, P, \{\cF\}_{t\ge 0})$ with $\B_0=\mathbf{0}$, $F_\x (\B)$ is a strong solution of \eqref{ISKE} with \eqref{:A1c} and \eqref{:A1d} starting at $\x$ for $P\circ \X_0^{-1}$-a.s. $\x$.
\end{dfn}

\subsection{Definition of \As[IFC]}\label{s:52}

In this subsection, we introduce \As[IFC] for our situation.
Let $\I$ be a finite subset of $\N$. Put $\I^c = \N \setminus \I$.
For $\y= (y^1, y^2, \dots) \in \S_{\rm hard}$, we put
$\y^{\I} = (y^j)_{j\in\I}$ and $\y^{\I^c} = (y^j)_{j\in\I^c}$
Let $(\X, \B)= ((X^j)_{j\in\N}, (B^j)_{j\in\N})$ be a weak solution of (\ref{ISKE}) starting at $\x=(x^j)_{j\in\N}$ defined on $(\Omega, \mathcal{F}, P, \{\mathcal{F}_t\})$.
We consider the SDE
\begin{align*} \nonumber %%%%%%%%%%%%%%%%
&dY_t^{\I,j} =  dB_t^j + 
b^{\I,j}_{\X} (t,\mathbf{Y}_t^{\I})dt
+\sum_{k\in \I\setminus \{j\}}(Y_t^{\I,j} - Y_t^{\I,k}) dL_t^{\I,jk}
 \\ \tag{SKE$_{\X}$($\I$)} \label{SKXI}
&\qquad\qquad\qquad+\sum_{k\in \I^c}^\infty (Y_t^{\I,j} - X_t^{k}) dL_t^{\I,jk},
\quad j\in \I,
\\
&Y_0^{\I,j} = X_0^j=x^j,
\quad j\in \I,
\end{align*}
where 
%$\displaystyle{b_{\X}^{\I,j}(t,\y)=b\bigg(y^j,\{y^k\}_{k\in\I \setminus \{j\}}+ \{X_t^k\}_{k\in \I^c}\bigg)}$ and 
$L_t^{\I,jk}$, $j\in\I$, $k\in \N$ are increasing functions satisfying
\begin{align*}
&L_t^{\I,jk}=\int_0^t \mathbf{1}(|Y_s^{\I,j}-Y_s^{\I,k}|=r)dL_s^{\I, jk},
\quad j,k \in \I,
\\
&L_t^{\I,jk}=\int_0^t \mathbf{1}(|Y_s^{\I,j}-X_s^k|=r)dL_s^{\I, jk},
\quad j\in \I, \ k\in \I^c.
\end{align*}

We denote by $\mathcal{C}^\I$ the completion of $\cB(W_{\bf 0}((\R^{d})^{\I})\times W((\R^d)^{\N}))$
with respect to $P_{\x}\circ (\B^\I, \X^{\I^c})^{-1}$.
Let $(\mathbf{v}, \w)\in W_{\bf 0}((\R^{d})^\I)\times W((\R^d)^{\N})$.
We %put $\cB_t=\sigma[\w_s : 0\le s \le t]$ and 
denote
by $\mathcal{C}_t^\I$ the completion of $\sigma[(\mathbf{v}_s, \w_s) : 0\le s \le t]$
with respect to $P_{\x}\circ (\B^\I, \X^{\I^c})^{-1}$.

%%%%%%%%%%%%%%%%%%%%%%%%%%%%%%%%%%%%%%%%%%%%%%
\begin{dfn}[strong solution for $(\X,\B)$ starting at $\x^\I$]\label{d:78b}
\quad % \\
$\mathbf{Y}^\I$ is called a strong solution of (\ref{SKXI}) for $(\X,\B)$ under   
$P_{\x}$ if $(\mathbf{Y}^\I, \B^{\I}, \X^{\I^c})$ satisfies 
(\ref{SKXI}) and there exists a $\mathcal{C}^\I$-measurable function 
$$
F_{\x}^\I :W_{\bf 0}((\R^d)^\I) \times W((\R^d)^{\I^c}) \to W((\R^d)^{\I})
$$
such that $F_{\bf x}^\I$ is $\mathcal{C}_t^\I / \cB_t(W(\R^d)^\I$-measurable for each $t$, and $F_{\x}^\I$ satisfies
$\mathbf{Y}^^I =F_{\x}^\I(\B^\I, \X^{\I^c})$, $P_\x$-a.s.
\end{dfn}
%%%%%%%%%%%%%%%%%%%%%%%%%%%%%%%%%%%%%%%%%%%%%%

%%%%%%%%%%%%%%%%%%%%%%%%%%%%%%%%%%%%%%%%%%%%%%
\begin{dfn}[a unique strong solution for $(\X,\B)$ starting at $\x_m$]\label{d:79b}
\quad \\
The SDE (\ref{SKXI}) is said to have a unique strong solution for ($\X, \B$) under $P_\x$ if there exists a strong solution $F_{\x}^\I$ such that for any solution $(\hat{\mathbf{Y}}^\I, \B^\I, \X^{\I^c})$ of (\ref{SKXI}) under $P_\x$,
%\begin{align}
$\hat{\mathbf{Y}}^\I = F_{\x}^\I(\B^\I, \X^{\I^c})\quad \mbox{ for $P_{\x}$- a.s.}$.
%\end{align}
\end{dfn}
%%%%%%%%%%%%%%%%%%%%%%%%%%%%%%%%%%%%%%%%%%%%%%
%\vskip 3mm
We can then give the definition of \As[IFC].

\vskip 3mm

\noindent \As[IFC] \quad For each finite subset $ \I \subset \N$, (\ref{SKXI}) has a unique strong solution
under $P_{\x}:= P(\cdot | \X_0 =\x)$ for $P \circ \X_0^{-1}$-a.s. $\x$.

%%%%%%%%%%%%%%%%%%%%%%%%%%%%%%%%%
\vskip 10mm
\noindent{\bf Acknowledgements.} \
This work was supported by JSPS KAKENHI Grant Numbers 
JP16H06338, JP19H01793, JP20K20K20855.
The author thanks Edanz \\
(https://jp.edanz.com/ac) for editing a draft of this manuscript.

%%%%%%%%%%%%%%%%%%%%%%%%%%%%%%%%

\end{document}